\documentclass[12pt]{amsart}

\usepackage{amssymb, amsmath, amsthm}
\usepackage{xcolor}
\definecolor{myred}{rgb}{0.2,0,0}
\definecolor{myblue}{rgb}{0,0,0.6}
\definecolor{mygreen}{rgb}{0,0.2,0}
\usepackage[colorlinks=true,linkcolor=myred,citecolor=myblue,urlcolor=mygreen]{hyperref}
\usepackage[centering]{geometry}

\newcommand{\les}{\leqslant}
\newcommand{\ges}{\geqslant}

\newtheorem{theorem}{Theorem}

\newtheorem{corollary}{Corollary}

\newtheorem{remark}{Remark}

\begin{document}
	
\title{Fractional parts of non-integer powers of primes} 
\author{Andrei Shubin}
\address{Department of Mathematics, Caltech, 1200 E. California Blvd., Pasadena, CA, 91125, USA}
\email{\href{mailto:ashubin@caltech.edu}{ashubin@caltech.edu}}
		
	\maketitle
		
	\begin{abstract}
		Let $\alpha > 0$ be any fixed non-integer, $I$ be any subinterval of $[0; 1)$. In the paper, we prove an analogue of Bombieri-Vinogradov theorem for the set of primes $p$ satisfying the condition $\{ p^{\alpha} \} \in I$. This strengthens the previous result of Gritsenko and Zinchenko.
	\end{abstract}

	\section{Introduction}
	\label{sec1}
	
	The problem concerning the distribution of primes satisfying the condition 
	$\{ p^{\alpha} \} <~\sigma$ for fixed numbers $\alpha$, $\sigma$ 
	such that $0 < \alpha, \sigma < 1$ was considered first by Vinogradov in~\cite{Vinogradov_1}.
	Using his method of exponential sums, he proved the formula
	\begin{equation} \label{Vinogradov_formula}
		\sum_{\substack{p \les X \\ \{p^{\alpha} \} < \sigma}} 1 = \sigma \pi(X) + R(X),
	\end{equation} where $\pi(X)$ is the prime-counting function and 
	$R(X) \ll X^{\vartheta(\alpha) + \varepsilon}$ with
	$$
		\vartheta(\alpha) = \max \biggl( \frac{4 + \alpha}{5}, 1 - \frac{2}{15} 
		\alpha  \biggr) = 
		\begin{cases}
		\displaystyle 1 - \frac{2}{15} \alpha, 
		&\text{if } \displaystyle 0 < \alpha \les 3/5;
		\\[2mm]
		\displaystyle\frac{4+\alpha}{5},
		&\text{if } \displaystyle 3/5 < \alpha < 1.
		\end{cases}
	$$ Linnik~\cite{Linnik} suggested another approach to this problem 
	based on zero density theorems for the Riemann zeta function. Using this 
	approach, Kaufman~\cite{Kaufman} proved the existence of infinite 
	number of primes $p$ such that $\{ \sqrt{p} \} < p^{-c}$ for any fixed $c$,
	$$
		0 < c < \frac{\sqrt{15}}{2 (8 + \sqrt{15})} = \frac{1}{6} - 0.00356\ldots .
	$$ This result was later extended to all $c < 1/4$ by Balog~\cite{Balog1}.
	
	The error term in~\eqref{Vinogradov_formula} was sharpened 
	by a number of authors (see~\cite{Balog1}--\cite{Ren}). 
	The best known result for $1/2 \les  \alpha < 1$ was obtained by Gritsenko~\cite{Gritsenko} via Linnik's approach:
	$$
		\vartheta(\alpha) = \max \biggl( \frac{4 + \alpha}{5}, 1 - \frac{2}{15} 
		\alpha  \biggr) = 
		\begin{cases}
		\displaystyle1 - \frac{\alpha}{2} + \bigl( \sqrt{3\alpha} - 1 \bigr)^2, 
		&\text{if } \displaystyle 1/2 \les  \alpha < 3/4;
		\\[2mm]
		\displaystyle\frac{1+\alpha}{2},
		&\text{if } \displaystyle 3/4 \les  \alpha < 1.
		\end{cases}
	$$ In~\cite{Harman} and~\cite{Gritsenko}, it was shown that one can take $\vartheta = 4/5$ 
	if $\alpha = 1/2$. For $0 < \alpha < 1/2$ the error term $R(X)$ was improved in~\cite{Ren}.
	
	An analogue of~\eqref{Vinogradov_formula} for large $\alpha$ was 
	obtained first by Vinogradov~\cite{Vinogradov_2}. He proved~\eqref{Vinogradov_formula} for all $\alpha > 6$ such that $|| \alpha || \ges 3^{-\alpha}$, 
	but with much worse power saving in the error term
	(here $|| x ||$ denotes the distance to nearest integer). This result was 
	extended to all non-integer $\alpha > 1$ and sharpened in~\cite{Leitmann}--\cite{Changa}. 
	
	Vinogradov gave an interesting interpretation of the subset
	of primes $p$ satisfying the restriction $\{ p^{\alpha} \} \in [c; d) \subset [0; 1)$: all 
	such primes lie in the intervals of the form 
	$$
		\bigl[ (k+c)^{1 / \alpha}; \ (k+d)^{1 / \alpha} \bigr).
	$$ Clearly, the length of such interval grows as $k \to +\infty$ for $\alpha < 1$ 
	and tends to zero for $\alpha > 1$. Thus, in the second 
	case the intervals are short and most of them do not contain even a single integer.
	This is the reason why the case $\alpha > 1$ seems to be more difficult.
	In particular, the Linnik approach does not work in this case.
	
	Many results in number theory are based on Bombieri-Vinogradov theorem,
	that is the estimate of the following type:
	\begin{equation} \label{BV}
		\sum_{q \les  Q} \max_{Y \les  X} \max_{(a,q)=1} \biggl| 
		\sum_{\substack{p \les  Y \\ p \equiv a\pmod{q}}} 1 - \frac{\pi(Y)}{\varphi(q)} 
		\biggr| \ll \frac{X}{(\log X)^A}.
	\end{equation} Here $A > 0$ is an arbitrary fixed constant, $Q = 
	X^{\theta - \varepsilon}$ with any fixed $\theta \les  1/2$. The exponent 
	$\theta$ is called usually as ``level of distribution''. 
	In the case of some
	special set $\mathbb{E}$ of integers, an analogue of Bombieri-Vinogradov 
	theorem has the form:
	\begin{equation} \label{analogue_BV}
	\sum_{q \les  Q} \max_{Y \les  X} \max_{(a,q)=1} \biggl| 
	\sum_{\substack{p 
			\les  Y, p \in \mathbb{E} \\ p \equiv a\pmod{q}}} 1 - 
	\frac{1}{\varphi(q)}  \sum_{\substack{p \les  Y \\ p \in \mathbb{E}}} 1 
	\biggr| \ll \frac{X}{(\log X)^A}.
	\end{equation} For the set $\mathbb{E} = \bigl\{ n \in \mathbb{N} : 
	\{ \sqrt{n} \} \in I \bigr\}$, where $I = [c;d) \subset [0;1)$, such an estimate was proved for the first time
	by Tolev~\cite{Tolev} with any $\theta \les  1/4$.
	Later, this result was 
	improved by Gritsenko and Zinchenko~\cite{Gritsenko-Zinchenko}, who 
	showed that~\eqref{analogue_BV} holds for all $\theta \les  1/3$ 
	for any set of the form
	$$
		\mathbb{E} = \bigl\{ n \in \mathbb{N} : \{ n^{\alpha} \} \in I\bigr\}, \qquad 1/2 \les  \alpha < 1.
	$$ In this work we obtain the similar estimate for all non-integer $\alpha > 0$. 
	
	The proof is based on the estimation of the exponential sums over primes of the form
	$$
		\sum_{\substack{X \les p < 2X \\ p\equiv a\pmod{q}}} e\bigl( hp^{\alpha} \bigr).
	$$ Here, as usual, $e(x) := e^{2 \pi i x}$. The main result is
	\begin{theorem} \label{thm1}
		Suppose that $\alpha > 0$ is fixed non-integer,
		$\theta, \varepsilon, D$ are fixed constants satisfying the 
		conditions $0 < \varepsilon < \theta < 1/3$, $\varepsilon < \alpha / 20$,
		$D > 1$, and suppose that 
		$1 \les  h \les  (\log X)^D$, 
		$2 < q \les  X^{\theta - \varepsilon}$, $1 \les a \les q-1$, $(a,q)=1$. Then the sum
		$$
			T = \sum_{\substack{X \les  p < 2X 
				\\ p \equiv a\pmod{q}}} e\bigl( h p^{\alpha}\bigr) 
		$$ satisfies the estimate
		$$
		T \ll \frac{X^{1 - \delta - \varepsilon^3 / (3 \alpha^2)}}{q},
		$$ where $0 < \delta \les  
		\varepsilon^3 / (50 \alpha^2)$, and the implied constant depends on $\alpha$, $\varepsilon$ and $D$.
	\end{theorem}

	As a corollary we deduce the formula for the proportion of primes $p \les X$ satisfying the restrictions $p \in \mathbb{E}, p \equiv a\pmod{q}$:

	\begin{corollary} \label{cor1}
		Let $\alpha > 0$ be a fixed non-integer, $\varepsilon > 0$ is arbitrary small number. Then for any $q \les X^{1/3 - \varepsilon}$, $a$, $(a,q)=1$, and any given subinterval $I \subset [0; 1)$ the following asymptotic formula holds true:
		$$
			\pi_{I} (X; q, a) := \sum_{\substack{p \les X \\ \{ p^{\alpha} \} \in I \\ p \equiv a\pmod{q}}} 1 = |I| \cdot \pi(X; q,a) + O \biggl( \frac{\pi(X; q,a)}{(\log X)^A} \biggr)
		$$ for any fixed $A > 0$.
	\end{corollary}

	Using \hyperref[thm1]{Theorem~1} we also obtain the analogue of Bombieri-Vinogradov theorem for primes in $\mathbb{E}$:
	
	\begin{theorem} \label{thm2} Suppose that $\alpha > 0$ is fixed non-integer and let
		$\mathbb{E}$ be the set of integers $n$ satisfying the condition
		$\{ n^{\alpha} \} \in I = [c; d) \subset [0; 1)$ for given $c$ and $d$. Further, let $\theta, \varepsilon$ and
		$A > 0$ be some fixed numbers such that $0 < \varepsilon < \theta < 1/3$, $\varepsilon < \alpha / 20$ and let
		$2 < Q \les  X^{\theta - \varepsilon}$. Then the inequality
		$$
		\sum_{q \les  Q} \max_{(a,q)=1} \biggl| \sum_{\substack{X \les  
				p < 2X \\ 
				p \equiv a\pmod{q} \\ p \in \mathbb{E}}} 1 - \frac{1}{\varphi(q)} 
		\sum_{\substack{X \les  p < 2X \\ p \in \mathbb{E}}} 1 \biggl| 
		\les  \frac{\kappa X}{(\log X)^A}
		$$ holds for any $X \ges X_0 (\alpha, \theta, \varepsilon)$ with
		some constant $\kappa > 0$ depending on $\alpha, \theta, \varepsilon$ 
		and $A$.
	\end{theorem}

	Recently, Bombieri-Vinogradov theorem was used for the proof of existence of the bounded gaps between primes (see~\cite{GPY}--\cite{Maynard}); for some generalization of this result, see~\cite{Maynard_2} and~\cite{Benatar}.

	In \hyperref[sec2]{Section 2}, we deduce \hyperref[thm2]{Theorem 2} from \hyperref[thm1]{Theorem 1} approximating the indicator function of the interval $[c; d)$ by ``Vinogradov's cup''. \hyperref[sec3]{Section 3} is devoted to the proof of \hyperref[thm1]{Theorem~1}. 
	In \hyperref[sec4]{Section 4}, we briefly discuss the application 
	of the main result to the bounded gaps between primes from $\mathbb{E}$. 
	
	Using Vaughan identity we split the initial sum
	\begin{equation} \label{exp_sum_general}
		\sum_{\substack{X \les  p < 2X \\ p \equiv a\pmod{q}}} e\bigl(h p^{\alpha}\bigr)
	\end{equation} to the sum of two types, namely,
	\begin{gather*}
	W_I = \sum_{M \les  m < 2M} \gamma_m \sum_{\substack{N \les  
			n < 2N \\ mn 
			\equiv a\pmod{q}}} f(n) e \bigl(h (mn)^{\alpha}\bigr), \\
	W_{II} = \sum_{M \les  m < 2M} \gamma_m \sum_{\substack{N \les  
			n < 2N \\ 
			mn \equiv a\pmod{q}}} \beta_n e \bigl(h (mn)^{\alpha}\bigr),
	\end{gather*} where the product $MN$ is close to $X$, $f(x)$ denotes some smooth 
	function and $\gamma_m$, $\beta_n$ denote some real coefficients.
	Then we obtain the appropriate upper bounds using the refined
	$k$-th derivative test for exponential sums, 
	which is due to Heath-Brown (see Theorem 1 in~\cite{Heath-Brown}).
	
	\begin{remark}
	In the forthcoming work we push the level of distribution $\theta$ in Theorems~\hyperref[thm1]{1} and \hyperref[thm2]{2} above the value $\theta = 1/3$ for small $\alpha$.
	\end{remark}

	\begin{remark}
	Recently in the paper ``A Bombieri-Vinogradov type result for exponential sums over primes'' Dimitrov~\cite{Dimitrov} claimed the level of distribution $1/2$ for exponential sums of the form~\eqref{exp_sum_general}, but with a tiny value of $h$, namely $h \ll X^{1/4 - \alpha}$ for $\alpha > 1$. Howewer, the case $h \ges 1$ is more interesting. 
	\end{remark}

	\section{Proof of \hyperref[thm2]{Theorem 2}}
	\label{sec2}
	
	Denote by $\chi(x)$ the indicator function of the interval $I = \bigl[ c; d 
	\bigr)$. Fixing some constant $B > 0$, we set $\Delta = (\log X)^{-B}$, 
	$r = \lfloor H \Delta \rfloor$ and
	$H = \Delta^{-1} \lceil \log_2 X \rceil$. Then there exists a 1-periodic
	function $\psi(x)$ (``Vinogradov's cup'', see~\cite[Ch.~1]{Karatsuba}) such that
	$\psi(x) = 1$ if $c+\Delta \les  x \les  d - \Delta$, $\psi(x) = 0$
	if $x \in [0; c] \cup [d; 1]$, $0 < \psi(x) < 1$ if 
	$x \in \bigl( c; c+\Delta \bigr) \cup \bigl(d - \Delta; d \bigr)$; moreover, $\psi(x)$ has the Fourier
	expansion of the form
	\begin{gather} \label{Fourier_expansion}
	\psi(x) = d - c - \Delta + \sum_{\substack{h = -\infty \\ h \ne 0}}^{
		+\infty} g(h) e(hx), \nonumber \\
	|g(h)| \les  \min \biggl( d - c - \Delta, \frac{1}{\pi |h|}, 
	\frac{1}{\pi |h|} \biggl( \frac{r}{\pi |h| \Delta} \biggr)^r \biggr).
	\end{gather} Therefore, setting
	$$
		\mathbb{E}_{\Delta} = \biggl\{ n \in \mathbb{N} : \{ n^{\alpha} \} \in \bigl( 
		c; c + \Delta \bigr) \cup \bigl( d - 
	\Delta; d \bigr)  \biggr\},
	$$ we obviously get
	$$
	\sum_{q \les  Q} \max_{(a,q) = 1} 
	\biggl| \sum_{\substack{X \les  p < 2X \\ n \equiv a\pmod{q}}} 
	\chi \bigl(\{ p^{\alpha} \}\bigr) - \frac{1}{\varphi(q)} \sum_{X \les  
		p < 2X} 
	\chi \bigl(\{ p^{\alpha} \}\bigr) \biggr| \les  S^{(1)} + S^{(2)} + 
	S^{(3)},
	$$ where
	\begin{gather*}
	S^{(1)} = \sum_{q \les  Q} \max_{(a,q) = 1} \biggl| 
	\sum_{\substack{X \les  p < 2X \\ p 
			\equiv a\pmod{q}}} \psi(p^{\alpha}) - \frac{1}{\varphi(q)} 
	\sum_{X \les  p < 2X} \psi(p^{\alpha}) \biggr|, \\
	S^{(2)} = \sum_{q \les  Q} \max_{(a,q)=1} \sum_{\substack{X \les  
			p < 2X \\ p 
			\equiv a\pmod{q} \\ p \in \mathbb{E}_{\Delta}}} 1, 
	\qquad S^{(3)} = \sum_{q \les  Q} \frac{1}{\varphi(q)} 
	\sum_{\substack{X \les  p < 2X 
			\\ p \in \mathbb{E}_{\Delta}}} 1.
	\end{gather*} Using~\eqref{Fourier_expansion}, we find
	\begin{multline} \label{S_1}
	S^{(1)} \les  \bigl( d - c - \Delta \bigr) \sum_{q \les  Q} 
	\max_{(a,q)=1} \biggl| \sum_{\substack{X \les  p < 2X \\ p \equiv 
			a\pmod{q}}} 1 - \frac{1}{\varphi(q)} \sum_{X \les  p < 2X} 1 \biggr| 
	+ \\ \biggl(\sum_{0 < |h| \les  H} + \sum_{|h| > H} \biggr) |g(h)|  
	\sum_{q \les  Q} \max_{(a,q)=1}
	\biggl| \sum_{\substack{X \les  p < 2X \\ p \equiv a\pmod{q}}} 
	e\bigl(h p^{\alpha}\bigr) 
	- \frac{1}{\varphi(q)} \sum_{X \les  p < 2X} e\bigl(h p^{\alpha}\bigr) \biggr|.
	\end{multline} By Bombieri-Vinogradov theorem, the first term in right-hand
	side in~\eqref{S_1} is estimated as
	$X (\log X)^{-C}$ for any fixed 
	$C > 0$. Trivial estimate of the inner sums over $p$ for $|h| > H$ together 
	with~\eqref{Fourier_expansion} yield:
	\begin{multline} \label{first_cup_big_h}
	\sum_{q \les  Q} \max_{(a,q)=1}  \sum_{|h| > H} |g(h)| \cdot \biggl| 
	\sum_{\substack{X \les  p < 2X \\ p \equiv a\pmod{q}}} e\bigl(h p^{\alpha}\bigr) - 
	\frac{1}{\varphi(q)} \sum_{X \les  p < 2X} e\bigl(h p^{\alpha}\bigr) \biggr| 
	\les  \\ 
    2X \sum_{q \les  Q} \sum_{|h| > H} \frac{1}{\pi |h|} \biggl( 
	\frac{r}{\pi 
		\Delta |h|}\biggr)^r \les  4X \sum_{q \les  Q} \frac{1}{\pi r} 
	\biggl( 
	\frac{r}{\pi \Delta (H-1)} \biggr)^r \les  X \sum_{q \les  Q} 
	\frac{1}{\pi} \biggl( 
	\frac{1}{2} \biggr)^r \les  Q.
	\end{multline} Next, the contribution coming from $0 < |h| \les  H$
	does not exceed
	$$
	\sum_{q \les  Q} \max_{(a,q)=1} \sum_{0 < |h| \les  H} |g(h)| 
	\biggl( \biggl| 
	\sum_{\substack{X \les  p < 2X \\ p \equiv a\pmod{q}}} e\bigl(h p^{\alpha}\bigr) 
	\biggr| + 
	\frac{1}{\varphi(q)} \biggl| \sum_{X \les  p < 2X} e\bigl(h p^{\alpha}\bigr) 
	\biggr| \biggr) 
	= S^{(4)} + S^{(5)}.
	$$ Using the estimates of the sum over primes $p$, $X \les  p < 2X$,
	given in~\cite{Vinogradov_1} (for $0 < \alpha < 1$) and~\cite{Changa}
	(for $\alpha > 1$), we get
	\begin{multline*}
	S^{(5)} \les  \sum_{q \les  Q} 
	\frac{1}{\varphi(q)} 
	\sum_{0 < |h| \les  H} \frac{1}{\pi 
		|h|} \biggl| \sum_{X \les  p < 2X} e\bigl(h p^{\alpha}\bigr) \biggr| 
	\ll_{\alpha, \varepsilon_1} \\
	X^{1 - \upsilon(\alpha) + \varepsilon_1} (\log Q) (\log H) \ll_{\alpha, \varepsilon_1} X^{1 - \upsilon(\alpha) + 2\varepsilon_1}
	\end{multline*} for arbitrarily small $\varepsilon_1 > 0$ and
	$$
	\upsilon(\alpha) = 
	\begin{cases}
	\displaystyle \alpha / 2,
	&\text{if } \displaystyle 0 < \alpha < 1;
	\\[2mm]
	\displaystyle\frac{6 \cdot 10^{-11}}{\alpha^2},
	&\text{if } \displaystyle \alpha > 1, \ \alpha \notin \mathbb{N}.
	\end{cases}
	$$ 
	Similarly, the estimate of \hyperref[thm1]{Theorem 1} yields:
	$$
	S^{(4)} \les  \sum_{0 < |h| \les  H} \frac{1}{\pi |h|} 
	\sum_{q \les  Q} 
	\max_{(a,q)=1} \biggl| 
	\sum_{\substack{X \les  p < 2X \\ p \equiv a\pmod{q}}} e\bigl(h 
	p^{\alpha}\bigr) \biggr| 
	\ll X^{1 - \varepsilon^3 / (3 \alpha^2)} \log H.
	$$ 
	
	Let $\Delta_1 = \Delta / 10$,  
	$r_1 = \bigl\lfloor H_1 \Delta_1 \bigr\rfloor$,  
	$H_1 = \Delta_1^{-1} \lceil \log_2 X \rceil$ and denote by $\psi_1(x)$
	and $\psi_2 (x)$ Vinogradov's cups such that $\psi_1 (x) = 1$ if
	$x \in (c; c + \Delta)$, $0 < \psi_1 (x) < 1$ if $x \in (c-\Delta_1; c) \cup
	(c+\Delta; c + \Delta + \Delta_1)$
	and $\psi_1 (x) = 0$ otherwise; $\psi_2 (x) = 1$
	if $x \in (d - \Delta; d)$, $0 < \psi_2 (x) < 1$ if 
	$x \in (d - \Delta - \Delta_1; d - \Delta) \cup (d; d + \Delta_1)$
	and $\psi_2 (x) = 0$ otherwise. Let us denote by $g_1 (h)$ and $g_2 (h)$
	its Fourier coefficients. Then
	\begin{multline} \label{R_1}
	S^{(2)} \les  \sum_{q \les  Q} \max_{(a,q)=1} \biggl| 
	\sum_{\substack{X \les  p < 2X \\ 
			p \equiv 
			a\pmod{q}}} \bigl(\psi_1 (p^{\alpha}) + \psi_2 (p^{\alpha})  \bigr) 
	\biggr| \les  \\ 2 (\Delta + \Delta_1) \sum_{q \les  Q}
	\max_{(a,q)=1} 
	\sum_{\substack{X \les  p < 2X \\ p \equiv a\pmod{q}}} 1 + 
	\sum_{h \ne 0} 	\biggl( |g_1 (h)| + |g_2 (h)| \biggr) 
	\sum_{q \les  Q} \max_{(a,q)=1} \biggl|  
	\sum_{\substack{X \les  p < 2X \\ p \equiv a\pmod{q}}} e\bigl(h p^{\alpha}\bigr) 
	\biggr|
	\end{multline} and, similarly,
	\begin{multline} \label{R_1'}
	S^{(3)} \les  2 (\Delta + \Delta_1 ) \sum_{q \les  Q} \frac{1}{\varphi(q)} 
	\sum_{X \les  p < 2X} 1 + \\
	\sum_{h \ne 0} \bigl( |g_1 (h)| + 
	|g_2 (h)| \bigr)
	\sum_{q \les  Q} \frac{1}{\varphi(q)} 
	\biggl| \sum_{X \les  p < 2X} e\bigl(h p^{\alpha}\bigr) \biggr|.
	\end{multline} Trivially, the first terms in the right hand side of~\eqref{R_1} and~\eqref{R_1'} do not exceed $2 \Delta X \log X$, and the second terms can be estimated similarly to $S^{(4)}$, $S^{(5)}$
	and the sum in the left hand side of~\eqref{first_cup_big_h}. To finish the proof, we choose $C = A, B = A + 1$.
	
	\begin{remark}
		In a similar way one deduces \hyperref[cor1]{Corollary~1} from \hyperref[thm1]{Theorem~1}.
	\end{remark}

	\section{Proof of \hyperref[thm1]{Theorem 1}}
	\label{sec3}
	
	Suppose that $1 \les  a < q \les  Q$, $(a,q)=1$, and consider the sum 
	$$
		W = W(Y) = \sum_{\substack{X \les  n < Y \\ n \equiv a\pmod{q}}} \Lambda(n) 
	e\bigl(hn^{\alpha}\bigr).
	$$ The application of 
	Vaughan identity with $V = X^{1/3}$ yields:
	$$
	W = -W_0 + W_1 - W_2 + W_3.
	$$ Here
	\begin{gather*}
	W_0 = \sum_{m \les  V^2} a_m \sum_{\substack{X/m \les  
			n < Y/m \\ mn \equiv 
			a\pmod{q}}} e\bigl( h (mn)^{\alpha}\bigr), \qquad 
	a_m = \sum_{\substack{uv = m \\ u,v \les  V}} \mu(u) \Lambda(v), \\
		W_1 = \sum_{\substack{n \les  V \\ n \equiv a\pmod{q}}} \Lambda(n) 
		e\bigl( h n^{\alpha}\bigr), \\
	W_2 = \sum_{V < m \les  Y V^{-1}} b_m \sum_{\substack{X/m \les  n 
			< Y/m, \ n > V 
			\\ mn \equiv a\pmod{q}}} \Lambda(n) e \bigl( h (mn)^{\alpha} \bigr), 
	\qquad b_m = \sum_{\substack{u|m \\ u \les  V}} \mu(u),	\\
	W_3 = \sum_{m \les  V} \mu(m) \sum_{\substack{X/m \les  n < Y/m \\ mn 
			\equiv 
			a\pmod{q}}} (\log n) e \bigl( h (mn)^{\alpha} \bigr).
	\end{gather*} Trivially, we have
	$$
		|a_m| \les  \sum_{v|m} \Lambda(v) = \log m, \qquad |b_m| \les  \tau(m),
	\qquad |W_1| \les  \sum_{n \les  V} \Lambda(n) \ll V.
	$$ Next, we have
	$$
	W_0 = \sum_{m \les  V} a_m \sum_{\substack{X/m \les  n < Y/m \\ 
			mn \equiv a\pmod{q}}} 
	e \bigl( h (mn)^{\alpha} \bigr) + \sum_{V < m \les  V^2} a_m 
	\sum_{\substack{X/m 
			\les  n < Y/m \\ mn \equiv a\pmod{q}}} e \bigl( h (mn)^{\alpha} 
	\bigr) = 
	W_4 + W_5.
	$$ Thus, we get type I sums $W_3$, $W_4$ and type II sums $W_2$, $W_5$. 
	Type I sums have the form
	$$
	\sum_{m \les  V} \gamma_m \sum_{\substack{X/m \les  n < Y/m \\ mn 
			\equiv 
			a\pmod{q}}} \beta_n e \bigl( h (mn)^{\alpha} \bigr),
	$$ where $\gamma_m = \mu(m)$, $\beta_n = \log n$ for $W_3$ and $\gamma_m = 
	a_m$, $\beta_n = 1$ for $W_4$; type II sums have the form
	$$
	\sum_{V < m \les  U} \gamma_m \sum_{\substack{Z/m \les  n < Y/m \\ 
			mn \equiv 
			a\pmod{q}}} \beta_n e \bigl( h (mn)^{\alpha} \bigr),
	$$ where $\gamma_m = b_m$, $\beta_n = \Lambda(n)$, $U = YV^{-1}$, $Z = 
	\max(Vm, X)$ for $W_2$ and $\gamma_m = a_m$, $\beta_n = 1$, $U = V^2$, 
	$Z = X$ for $W_5$.

	\subsection*{The estimation of type I sums}
	We split the range of summation 
	$1 \les  m \les  V$ to the dyadic intervals 
	$M < m \les  M_1$, $M_1 = \min (2M, V)$. Then the initial sum splits 
	into $\ll \log X$ sums of the form
	$$
	W(M) = \sum_{M < m \les  M_1} \gamma_m \sum_{\substack{X/m \les  
			n < Y/m \\ mn 
			\equiv a\pmod{q}}} \beta_n e \bigl( h (mn)^{\alpha} \bigr).	
	$$ By partial summation, we get
	\begin{multline} \label{Type_I_after_Abel}
	\bigl|W(M)\bigr| \les  2 ||\beta||_{\infty} \sum_{\substack{M < m 
			\les  M_1 \\ 
			(m,q) = 1}	} \bigl|
	\gamma_m\bigr| \cdot \biggl| \sum_{\substack{X/m \les  n < Y_1/m \\ mn 
			\equiv a\pmod{q}}} e \bigl( h (mn)^{\alpha} \bigr) \biggr| \les  
	\\ 2 || \gamma ||_{\infty} || \beta ||_{\infty} \sum_{\substack{M 
			< m \les  M_1 \\ (m,q)=1}} \biggl| \sum_{\substack{X/m 
			\les  n < 
			Y_1/m \\ mn \equiv a\pmod{q}}} e \bigl( h (mn)^{\alpha} \bigr) 
	\biggr|,
	\end{multline} where $Y_1 \in (X; Y]$ and $|| \omega ||_{\infty} = \max_{n 
		\les  2X} |\omega_n|$. Next, we fix $m \in 
	(M; M_1]$ with $(m,q) = 1$ and define $l \equiv am^{\ast} \pmod{q}$, $1 
	\les  l \les  q-1$. Setting $n = qr + l$ we obtain
	$$
	\frac{X}{mq} \les  r + \xi < \frac{Y_1}{mq}, \qquad \xi = \frac{l}{q}.
	$$ The inner sum over $n$ in~\eqref{Type_I_after_Abel} takes the form
	\begin{equation} \label{Type_I_after_sub}
	\sum_{R_1 - \xi \les  r < R_2 - \xi} e \bigl( h (mq)^{\alpha} 
	(r + \xi)^{\alpha} 
	\bigr), 
	\end{equation} where $R_1 = X / mq, R_2 = Y_1 / mq \les  2R_1$. 
	
	To estimate aforementioned sum, we need the following assertion 
	(Theorem 1,~\cite{Heath-Brown}): suppose that $k \ges 3$ and let $f \in C^k 
	\bigl( [N; 2N] \bigr)$ satisfies the conditions
	$0 < \lambda_k \les  f^{(k)}(x) \les  L \lambda_k$ for any 
	$x \in [N; 2N]$. Then, for any fixed $\delta > 0$, we have
	\begin{multline} \label{Heath-Brown}
	\sum_{N \les  n < 2N} e\bigl(f(n)\bigr) \ll_{L, k, \delta} 
	\\ N^{1 + 
		\delta} \bigl(  \lambda_k^{1 / (k(k-1))} + N^{-1 / (k(k-1))} + 
	N^{-2 / (k(k-1))} \lambda_k^{-2 / (k^2 (k-1))} \bigr).
	\end{multline}
	
	Consider the 
	function $f_I (x) = h (mq)^{\alpha} (x + \xi)^{\alpha}$. Then, for 
	$R_1 - \xi \les  x < R_2 - \xi$,
	$$
	f_I^{(k)} (x) = (\alpha)_k h (mq)^{\alpha} (x + \xi)^{\alpha - k} \asymp \frac{h 
		(mq)^k}{X^{k-\alpha}} = \lambda_k,
	$$ where $(\alpha)_k = \prod_{i=1}^k (\alpha - i + 1)$ is the Pochhammer 
	symbol. Applying~\eqref{Heath-Brown} to~\eqref{Type_I_after_sub}, we get
	\begin{multline*}
	\biggl| \sum_{R_1 - \xi \les  r < R_2 - \xi} e\bigl( f_I (r) \bigr) 
	\biggr| \ll_{k, 
		\delta} \biggl( \frac{X}{mq} \biggr)^{1 + \delta} \cdot   
	\biggl[ \biggl( \frac{h (mq)^k}{X^{k-\alpha}} \biggr)^{1 / (k(k-1))} 
	+ \biggl( \frac{mq}{X} \biggr)^{1/(k(k-1))} + \\ \biggl( \frac{mq}{X} 
	\biggr)^{2/(k(k-1))} \biggl( \frac{X^{k-\alpha}}{h(mq)^k} \biggr)^{2 / (k^2 
		(k-1))} \biggr] \ll
	h^{1/(k(k-1))} \biggl( \frac{X}{mq} \biggr)^{\delta} \biggl( T_1 + 
	T_2 + T_3 \biggr),
	\end{multline*} where
	$$
	T_1 = \frac{ X^{1 - (k-\alpha) / (k(k-1))}}{ (mq)^{1 - 1/(k-1)}}, \qquad
	T_2 = \biggl(\frac{ X}{mq} \biggr)^{1 - 1/(k(k-1))}, \qquad T_3 = \frac{X^{1 
			- 2\alpha / (k^2 (k-1))}}{mq}.
	$$ Since $||\gamma||_{\infty} || \beta ||_{\infty} \les  \log 2X$, 
	we get
	$$
	\bigl|W(M)\bigr| \ll h^{1 / (k(k-1))} (\log X) \biggl( \frac{X}{q} 
	\biggr)^{\delta} \sum_{M < m \les  M_1} m^{-\delta} 
	( T_1 + T_2 + T_3). 
	$$ Now we estimate the contribution from $T_1, T_2, T_3$ to the sum 
	over all values of $M$. The contribution from $T_1$ doesn't exceed
	\begin{multline} \label{Type_I_T_1}
	h^{1/(k(k-1))} (\log X) \biggl( \frac{X}{q} \biggr)^{\delta} \frac{X^{1 
			- (k-\alpha)/(k(k-1))}}{q^{1-1/(k-1)}} {\sum_{M < V}}' \sum_{M < m 
		\les  M_1} m^{-1 + 1/(k-1)} \ll \\ \biggl( \frac{X}{q} \biggr)^{2
		\delta} \frac{X^{1 - (k-\alpha)/(k(k-1))}}{q^{1-1/(k-1)}} 
	V^{1/(k-1)}.
	\end{multline} The contribution from $T_2$ is less than
	\begin{multline} \label{Type_I_T_2}
	h^{1/k(k-1)} (\log X) \biggl( \frac{X}{q} \biggr)^{1 - 1/(k(k-1)) + 
		\delta} {\sum_{M < V}}' \sum_{M < m \les  M_1} m^{-1 + 
		1/(k(k-1))} \ll \\ \biggl( \frac{X}{q} \biggr)^{1 - 1/(k(k-1)) 
		+ 2 \delta} V^{1/(k(k-1))}.
	\end{multline} Finally, the contribution from $T_3$ is bounded by
	\begin{multline}\label{Type_I_T_3}
	h^{1/(k(k-1))} (\log X) \biggl( \frac{X}{q} \biggr)^{\delta} 
	\frac{X^{1 - 2\alpha / (k^2 (k-1))}}{q} {\sum_{M < V}}' 
	\sum_{M < m \les  M_1} 
	\frac{1}{m} \ll \\  \biggl( \frac{X}{q} \biggr)^{2 \delta} 
	\frac{X^{1 - 2\alpha / (k^2 (k-1))}}{q} \log V. 
	\end{multline}

	\subsection*{The estimation of type II sums}
	By definition of $U$, $V^2 
	\les  U \les  Y V^{-1} \les  2 X V^{-1} \les  2V^2$. We split 
	$W_2$, $W_5$
	into $\ll \log X$ sums of the type $W(M)$. Cauchy's inequality yields:
	$$
	\bigl|W(M)\bigr|^2 \les  \biggl( \sum_{M < m \les  M_1} |\gamma_m|^2 
	\biggr) 
	\biggl( \sum_{M < m \les  M_1} \biggl| \sum_{\substack{Z/m \les  n < 
			Y/m \\ mn 
			\equiv a\pmod{q}}} \beta_n e \bigl( h (mn)^{\alpha} \bigr) \biggr|^2 
	\biggr).		
	$$ Next, by Mardzhanishvili's inequality~\cite{Mardzhanishvili} we get
	\begin{equation} \label{Mardz}
	\bigl|W(M)\bigr|^2 \ll M (\log X)^{2 + \kappa}
	\biggl( \sum_{M < m \les  M_1} \biggl| \sum_{\substack{Z/m \les  n < 
			Y/m \\ mn 
			\equiv a\pmod{q}}} \beta_n e \bigl( h (mn)^{\alpha} \bigr) \biggr|^2 
	\biggr), 
	\end{equation} where $\kappa = 1$ for $W_2$ and $\kappa = 0$ for $W_5$. Now we rewrite 
	the sum over $m$ as follows:
	\begin{multline*}
	\sum_{M < m \les  M_1} \sum_{\substack{Z/m \les  n_1, n_2 < 
			Y/m \\ mn_i \equiv 
			a\pmod{q}, i=1,2}} \beta_{n_1} \beta_{n_2} e \bigl( h 
	m^{\alpha} (n_1^{\alpha} 
	- n_2^{\alpha}) \bigr) = \\ \sum_{M < m \les  M_1} \sum_{\substack{Z/m 
			\les  n < Y/m \\ mn \equiv a\pmod{q}}} {\beta^2_n} + 2 \text{Re} 
	(S(M)),
	\end{multline*} where
	$$
	S(M) = \sum_{M < m \les  M_1} \sum_{\substack{Z/m \les  n_1 < n_2 < 
			Y/m \\ mn_i 
			\equiv a\pmod{q}, i=1,2}} \beta_{n_1} \beta_{n_2} e \bigl( h m^{\alpha} 
	(n_1^{\alpha} - n_2^{\alpha}) \bigr).
	$$ The diagonal term doesn't exceed
	\begin{multline} \label{diagonal}
	\sum_{M < m \les  M_1} \sum_{\substack{Z/m \les  n < Y/m \\ 
			mn \equiv a\pmod{q}}} 
	\beta_n^2 \ll \sum_{M < m \les  M_1} (\log X)^{2 \kappa} \biggl( \frac{X}{mq} + 
	1 \biggr) \ll \\ (\log X)^{2 \kappa} \biggl( \frac{X}{q} + M \biggr).
	\end{multline} Setting 
	$m = qr + l$, we get
	$$
	\frac{M}{q} - \eta < r \les  \frac{M_1}{q} - \eta, \qquad \eta = \frac{l}{q}
	$$ for given $l$, $(l,q) = 1$. Hence,
	$$
	S(M) = \sum_{\substack{l=1 \\ (l,q)=1}}^q \sum_{\frac{M}{q} - \eta < r \les 
		\frac{M_1}{q} - \eta} \sum_{\substack{\frac{Z}{qr+l} \les  n_2 < n_1 < 
			\frac{Y}{qr+l} \\ n_1, n_2 \equiv e\pmod{q}}} \beta_{n_1} \beta_{n_2}  
	e \bigl( h (n_1^{\alpha} - n_2^{\alpha}) q^{\alpha} (r + \eta)^{\alpha} 
	\bigr),
	$$ where $e = a l^{\ast} \pmod{q}$. Next, we change the order 
	of summation. If $Z = X$ then
	$$
	\frac{X}{qr + l} \les  n_2 < n_1 < \frac{Y}{qr + l},
	$$ so $X/M_1 \les  n_2 < n_1 < Y/M$ and for fixed $n_1$, $n_2$ we get
	$$
	\frac{X}{qn_2} - \eta \les  r < \frac{Y}{qn_1} - \eta.
	$$ By definition, $Z = \max(Vm, X) = \max \bigl( V(qr+l), X \bigr)$ so hence
	$$
	\max \biggl( V , \frac{X}{qr+l} \biggr) \les  n_2 < n_1 < \frac{Y}{qr+l}.
	$$ Since
	$$
	\max \biggl( V, \frac{X}{qr+l} \biggr) = 
	\begin{cases}
	\displaystyle\frac{X}{qr+l},
	&\text{if }\displaystyle r \les  \frac{XV^{-1}}{q} - \eta;
	\\[2mm]
	V,
	&\text{if }\displaystyle r > \frac{XV^{-1}}{q} - \eta.
	\end{cases}
	$$ then we estimate $S(M)$ as follows:
	\begin{multline*}
	S(M) = \biggl\{ \sum_{\substack{\frac{M}{q} - \eta < r \les  \frac{M_1}{q} - 
			\eta \\ 
			r \les  XV^{-1} / q - \eta}} \sum_{\substack{\frac{X}{qr+l} 
			\les  n_2 < 
			n_1 < \frac{Y}{qr + l} \\ n_1, n_2 \equiv e\pmod{q}}} + \sum_{
		\substack{\frac{M}{q} - \eta < r \les  \frac{M_1}{q} - 
			\eta \\ r > XV^{-1} 
			/ q - \eta}} \sum_{\substack{V \les  n_2 < n_1 < \frac{Y}{qr + l}
			\\ 
			n_1, n_2 \equiv e\pmod{q}}} \biggr\} \ldots = \\ \biggl\{ \sum_{
		\substack{X/M_1 \les  n_2 < n_1 < Y/M \\ n_1, n_2 \equiv e\pmod{q}}} 
	\sum_{R^{(1)} - \eta < r \les  R^{(2)} - \eta} + \sum_{\substack{V 
			\les  n_2 
			< n_1 < Y/M \\ n_1, n_2 \equiv e\pmod{q}}} \sum_{R^{(3)} - \eta < 
		r \les  R^{(4)} - \eta}\biggr\} \ldots,
	\end{multline*} where
	\begin{gather*}
	R^{(1)} = \max \biggl( \frac{M}{q}, \frac{X}{qn_2} \biggr), \qquad R^{(2)} 
	= \min \biggl( \frac{M_1}{q}, \frac{Y}{qn_1}, \frac{XV^{-1}}{q} \biggr), \\
	R^{(3)} = \max \biggl( \frac{M}{q}, \frac{X V^{-1}}{q} \biggr), \qquad R^{(4)} 
	= \min \biggl( \frac{M_1}{q}, \frac{Y}{qn_1} \biggr).
	\end{gather*} Therefore,
	$$
	\bigl|S(M)\bigr| \les  \sum_{\substack{l=1 \\ (l,q)=1}}^q 
	\sum_{\substack{X / 
			M_1 \les  n_2 < n_1 < Y/M \\ n_1, n_2 \equiv e\pmod{q} }} 
	|\beta_{n_1}| 
	|\beta_{n_2}| \biggl| \sum_{R_1 - \eta < r \les  R_2 - \eta} 
	e\bigl( f_{II}(r)\bigr) 
	\biggr|,
	$$ where $(R_1, R_2)$ denotes the pair $(R^{(1)}, R^{(2)})$, $(R^{(3)}, 
	R^{(4)})$ 
	that corresponds to the maximum absolute value of the sum over $r$, 
	$f_{II} (x) = 
	h(n_1^{\alpha} - n_2^{\alpha}) q^{\alpha} (x + \eta)^{\alpha}$. Using the 
	conditions $n_2 < n_1$, 
	$n_1 \equiv n_2 \equiv e \pmod{q}$ we write $n_1 = 
	n_2 + qs$ with $s \ges 1$. On the other hand, $n_1 < Y / M$ implies $n_2 + 
	qs < Y/M$. Hence, $s < Y / (Mq) = t$ and therefore
	$$
	\bigl|S(M)\bigr| \ll \sum_{\substack{l=1 \\ (l,q)=1}}^q \sum_{1 \les  s < t} 
	\sum_{\substack{X/M_
			1 \les  n < Y/M \\ n \equiv e\pmod{q}}} |\beta_n| |\beta_{n+qs}| 
	\biggl| \sum_{R_1 \les  r < R_2} e \bigl( f_{II}(r) \bigr) \biggr|.
	$$ Obviously,
	$$
	f_{II}^{(k)} (x) = \frac{(\alpha)_k h (n_1^{\alpha} - n_2^{\alpha}) q^{
			\alpha}}{(x + \eta)^{k-\alpha}} = \frac{(\alpha)_k D_1}{(x + \eta)^{
			k-\alpha}},
	$$ where
	$$
	D_1 = h (n_1^{\alpha} - n_2^{\alpha}) q^{\alpha}, 
	\qquad \text{so we have} \quad
	\bigl|f_{II}^{(k)} (x)\bigr| \asymp \frac{D_1}{R_1^{k-\alpha}} \asymp D_1 \biggl( 
	\frac{q}{M} 
	\biggr)^{k-\alpha}.
	$$ Next, by Lagrange's mean value theorem,
	$$
	D_1 = hq^{\alpha} \bigl( (n + qs)^{\alpha} - n^{\alpha} \bigr) = hq^{\alpha} 
	\cdot \alpha (n + qs \theta')^{\alpha - 1} \cdot qs \asymp hs q^{\alpha + 1} 
	\biggl( \frac{X}{M} \biggr)^{\alpha - 1}, \qquad |\theta'| \les  1,
	$$ and hence
	$$
	\bigl|f_{II}^{(k)}(x)\bigr| \asymp \frac{hsq^2}{X^{1 - \alpha}} \biggl( \frac{q}{M} 
	\biggr)^{k-1} = \lambda_k.
	$$ Put $1 - \alpha = \nu$. By~\eqref{Heath-Brown}, we get
	\begin{multline*}
	\sum_{R_1 < r \les  R_2} e \bigl( f_{II}(r) \bigr) \ll_{k, \delta} 
	\biggl( \frac{M}{q} \biggr)^{1 + \delta} \biggl\{ \biggl( 
	\frac{hsq^2}{X^{\nu}} \biggr)^{1 / (k(k-1))}  
	\biggl( \frac{q}{M} \biggr)^{1/k} + \biggl( \frac{M}{q} \biggr)^{
		-1 / (k(k-1))} + \\ \biggl( \frac{M}{q} \biggr)^{-2 / (k(k-1))} 
	\biggl( \frac{hsq^2}{X^{\nu}} \biggr)^{-2 / (k^2 (k-1))}  
	\biggl( \frac{q}{M} \biggr)^{-2/k^2} \biggr\}. 
	\end{multline*} The factor $|\beta_n| \cdot |\beta_{n+qs}|$ is bounded 
	from above by $(X/q)^{\delta}$. 
	The summation over
	$n \equiv e\pmod{q}$ for $X/M_1 < n \les  Y/M$ contributes the 
	factor of at most $X / Mq$. Thus,
	\begin{multline*}
	S(M) \ll \biggl( \frac{X}{q} \biggr)^{\delta} \sum_{\substack{l=1 \\ 
			(l,q)=1}}^q \frac{X}{Mq} \biggl( \frac{M}{q} \biggr)^{1+\delta} 
	\biggl\{ \biggl( \frac{hq^2}{X^{\nu}} \biggr)^{1 / (k(k-1))} 
	\biggl( \frac{q}{M} \biggr)^{1/k} \sum_{1 
		\les  s < t} s^{1/(k(k-1))} + \\ \biggl( \frac{M}{q} \biggr)^{
		-1/(k(k-1))} \sum_{1 \les  s < t} 1 + \biggl( \frac{hq^2}{X^{\nu}} 
	\biggr)^{-2 / (k^2 (k-1))} \biggl( \frac{M}{q} 
	\biggr)^{-2 / (k^2 (k-1))} \sum_{1 \les  s < t} s^{-2 / (k^2 (k-1))} 
	\biggr\}.
	\end{multline*} The inequalities $M < X$ and $t \les  2X / Mq$ imply:
	\begin{multline} \label{T_4_T_5_T_6}
	S(M) \ll \biggl( \frac{X}{q} \biggr)^{2 \delta} \cdot \frac{X}{q} 
	\biggl\{ \biggl( 
	\frac{hq^2}{X^{\nu}} \biggr)^{1/(k(k-1))} \biggl( 
	\frac{q}{M} \biggr)^{1/k} \biggl( \frac{2X}{Mq} \biggr)^{1 + 
		1/(k(k-1))} + \\ \biggl( \frac{M}{q} \biggr)^{-1 / (k(k-1))} 
	\frac{2X}{Mq} + \biggl( \frac{hq^2}{X^{\nu}} \biggr)^{-2 / (k^2 (k-1))} 
	\biggl( \frac{M}{q} \biggr)^{
		-2 / (k^2 (k-1))} \biggl( \frac{2X}{Mq} \biggr)^{1 - 2 /
		(k^2 (k-1))} \biggr\} \ll \\ \biggl( \frac{X}{q} \biggr)^{2 \delta} 
	\frac{X^2}{Mq^2} (T_4 + T_5 + T_6),
	\end{multline} where
	\begin{gather*}
	T_4 = (2h)^{1 / (k(k-1))} X^{(1-\nu) / (k(k-1))} \biggl(\frac{q}{M}
	\biggr)^{1 / (k-1)} 
	= \biggl( \frac{2h 
		q^k  X^{\alpha}}{M^k} \biggr)^{
		1 / (k(k-1))}, \\
	T_5 = \biggl( \frac{q}{M} \biggr)^{1 / (k(k-1))}, \qquad
	T_6 = (2h)^{-2 / (k^2 (k-1))} X^{(2\nu - 2) / (k^2 (k-1))} 
	= (2 h X^{\alpha})^{-2 / (k^2 (k-1))}.  
	\end{gather*} Thus, the contribution from~\eqref{T_4_T_5_T_6} to $|W(M)|^2$ doesn't exceed
	$$
	(\log X)^{2 + \kappa} \biggl( \frac{X}{q} \biggr)^{2 + 2 \delta} 
	\biggl\{  \biggl( \frac{2h 
		q^k  X^{\alpha}}{M^k} \biggr)^{
		1 / (k(k-1))} + 
	\biggl( \frac{q}{M} \biggr)^{1 / (k(k-1))} + \biggl( \frac{1}{2h X^{\alpha}} 
	\biggr)^{2 / (k^2 (k-1))} \biggr\},
	$$ whence, combining with~\eqref{Mardz} and~\eqref{diagonal} we get
	\begin{multline*}
	\bigl|W(M)\bigr| \ll_k (\log X)^{1 + 3\kappa / 2} \biggl( \frac{X}{q} \biggr)^{
		\delta} \biggl\{ M + \biggl(\frac{MX}{q} \biggr)^{1/2} + \\
	\frac{X}{q} \biggl( \biggl( \frac{h 
		q^k  X^{\alpha}}{M^k} \biggr)^{
		1 / (2k(k-1))} + \biggl( \frac{q}{M} \biggr)^{1 / (2k (k-1))} + 
	\biggl( \frac{1}{h X^{\alpha}} \biggr)^{1 / (k^2 (k-1))} \biggr) \biggr\}.
	\end{multline*} The summation over all $M$ in the range $V 
	\les  M < 2V^2$ leads to the estimate
	\begin{multline} \label{Type_II_total}
	W \ll (\log X)^{1 + 3 \kappa / 2} \biggl( \frac{X}{q} \biggr)^{\delta} 
	\biggl\{ 
	V^2 + \frac{V \sqrt{X}}{\sqrt{q}} + \\ \frac{X}{q} \biggl( \biggl( \frac{h 
		q^k  X^{\alpha}}{V^k} \biggr)^{
		1 / (2k(k-1))} + \biggl( 
	\frac{q}{V} \biggr)^{1 / (2k(k-1))} + \biggl( \frac{1}{h X^{\alpha}} 
	\biggr)^{1 / (k^2 (k-1))} \biggr) \biggr\}.
	\end{multline}

	\subsection*{Final bound}
	From~\eqref{Type_I_T_1}, \eqref{Type_I_T_2}, 
	\eqref{Type_I_T_3}, \eqref{Type_II_total} we conclude that
	\begin{multline*}
	W \ll \biggl( \frac{X}{q} \biggr)^{1 + 2\delta} \biggl\{
	\biggl(\frac{Vq}{X} \biggr)^{1 / (k-1)} X^{\alpha / (k(k-1))} +
	\biggl( \frac{Vq}{X} \biggr)^{1 / (k(k-1))} + X^{-2\alpha / (k^2 (k-1))}
	\log V  
	+ \\ \frac{V^2 q}{X} + \frac{V \sqrt{q}}{\sqrt{X}} + 
	\biggl( \frac{h q^k X^{\alpha}}{V^k} \biggr)^{
		1 / (2k(k-1))} + \biggl( \frac{q}{V} \biggr)^{1 / (2k (k-1))} + 
	\biggl(\frac{1}{hX^{\alpha}}\biggr)^{1 / (k^2 (k-1))}  \biggr\}.
	\end{multline*} We estimate the factors 
	$( X/q )^{\delta}$ and $h^{1 / (2k(k-1))}$ by $X^{\delta}$. 
	Thus,
	\begin{multline*}
	W \ll X^{1 + 3\delta}  \biggl\{ 
	\biggl(\frac{Vq}{X^{1-\alpha / k}} \biggr)^{1 / (k-1)} 
	+ \biggl(\frac{Vq}{X}\biggr)^{1 / (k(k-1))} + 
	X^{-2\alpha / (k^2 (k-1))}  + 
	\frac{ V^2 q}{X} + 
	\frac{V\sqrt{q}}{\sqrt{X}} + \\ 
	\biggl( 
	\frac{q^k X^{\alpha}}{V^k} \biggr)^{1 / (2k(k-1))} + \biggl( 
	\frac{q}{V} \biggr)^{1 / (2k(k-1))} + \biggl( \frac{1}{X^{\alpha}} 
	\biggr)^{1 / (k^2 (k-1))} \biggr\} \ll X^{1 + 3\delta} \sum_{i=1}^8 
	\Delta_i,
	\end{multline*} where
	\begin{gather*}
	\Delta_1 \les  X^{(3\alpha - 2k) / (3k(k-1))} q^{1 / (k-1)} 
	\ll X^{(3 \alpha - k - 3 \varepsilon k) / (3k (k-1))}, \qquad
	\Delta_2 = \biggl(\frac{q}{X^{2/3}}\biggr)^{1 / (k(k-1))}, \\
	\Delta_2 \ll X^{-1 / (3k(k-1))}, \ \
	\Delta_3 = X^{-2 \alpha / (k^2 (k-1))}, \ \
	\Delta_4 = \frac{V^2 q}{X} \les  X^{-\varepsilon}, \ \
	\Delta_5 = \frac{V \sqrt{q}}{\sqrt{X}} \les  X^{-\varepsilon / 2}, \\
	\Delta_6 \les  \bigl( X^{\alpha / k - 
		\varepsilon} \bigr)^{1 / (2(k-1))}, \qquad 
	\Delta_7 \les  X^{-\varepsilon / (2k(k-1))}, \qquad \Delta_8 \les  X^{
		- \alpha / (k^2 (k-1))}, \\
	\max_{1 \les  i \les  8} \Delta_i \les  X^{ - 
		2\varepsilon^3 / (5 \alpha^2)}
	\end{gather*} if $k = \lfloor 1.1 \cdot \alpha / \varepsilon \rfloor + 
	1$. Finally, choosing
	$\delta \les  \varepsilon^3 / (50 \alpha^2)$ and applying partial summation,
	we lead to the desired bound.

	\section{Application to bounded gaps}
	\label{sec4}
	
	In this section, we prove the existence of infinite number of bounded gaps between successsive primes from $\mathbb{E}$. We follow the well-known technique of Maynard~\cite{Maynard} and Tao with modified Selberg sieve and 
	\hyperref[thm2]{Theorem 2} in place of the Bombieri-Vinogradov theorem. We consider only the case $0 < \alpha < 1$, which is more simple. With a little more effort, one can deduce a similar result for any non-integer $\alpha > 1$.
	
	\begin{theorem}\label{thm3}
		Let $\mathbb{E} = \bigl\{ n \in \mathbb{N} : \{ n^{\alpha} \} \in [c; d) \subset [0; 1) 
		\bigr\}$ for given $c$ and $d$, $0 < \alpha < 1$, $q_1, q_2, \ldots, q_n, \ldots$ be all primes 
		from $\mathbb{E}$ indexed in ascending order, and suppose that
		$m \ges 1$ is a fixed integer. 
		Then
		$$
			\liminf_{n \to +\infty} (q_{n+m} - q_n) \les  9~700 m^3 e^{6m}.
		$$
	\end{theorem}
	
	The set $\{h_1, \ldots, h_k \}$ of integers is called an \textit{admissible set}
	if for any prime $p$ there 
	is an $a$ with $h_j\not\equiv a\pmod{p}$ for any $1 \les  j \les  k$.
	Consider the sum
	$$
	S_{\alpha} = \sum_{\substack{X \les  n < 2X \\ n \in \mathbb{E} \\ n+h_j \in 
			\mathbb{E} \ \forall j}} \biggl( \sum_{j=1}^k \chi_{ 
		\mathbb{P}} (n+h_j) - \rho \biggr) \omega_n = S_{2,\alpha} - \rho 
	S_{1,\alpha},
	$$ where $\omega_n, \rho > 0$, $\chi_{\mathbb{P}}$ is the 
	characteristic function of primes. Then 
	\hyperref[thm3]{Theorem 3} clearly 
	follows from the inequality $S_{2,\alpha} - \rho S_{1,\alpha} > 0$ with 
	$\rho = m$. Indeed, in this case the inner 
	sum over $1 \les  j \les  k$ has at least $m+1$ positive terms for 
	some~$n$. Hence, there are at least $m+1$ primes from $\mathbb{E}$ between
	$n$ and $n+h_k$ and
	\begin{equation} \label{maximum}
	\liminf_{n \to +\infty} (q_{n+m} - q_n) \les  
	\max_{1 \les  i < j \les  k} 
	|h_j - h_i| .
	\end{equation} Thus, the problem reduces to choosing the appropriate weights 
	$\omega_n$ 
	maximizing the ratio $S_{2,\alpha} / S_{1,\alpha}$. We follow the choise made 
	in~\cite{Maynard}:
	$$
	\omega_n = \biggl( \sum_{d_j | (n+h_j)} \lambda_{d_1, \ldots, d_k} \biggr)^2,
	$$ where the sum is taken over all tuples of divisors $(d_1, \ldots, d_k)$ 
	and
	$$
	\lambda_{d_1, \ldots, d_k} = \biggl( \prod_{j=1}^k \mu(d_j) d_j \biggr)  
	\sum_{\substack{r_1, \ldots, r_k \\ d_j | r_j \ \forall j \\ (r_j, W) = 1 
			\ \forall j}} \frac{\mu^2 
		(r_1 \ldots r_k)}{\varphi(r_1)\ldots \varphi(r_k)} F \biggl( \frac{\log 
		r_1}{\log R}, \ldots, \frac{\log r_k}{\log R} \biggr).
	$$ Here $W$ is the product of all primes $\les  \log \log \log X$ 
	(and so $W 
	\les  (\log \log X)^2$ for large $X$), $R = X^{1/6 - \delta_1}$ for 
	some small fixed $\delta_1 > 0$ and $F(x_1, \ldots, x_k)$ is a fixed 
	piecewise continuous function supported on the set 
	$$
		\biggl\{ (x_1, \ldots, 
	x_k) \in [0,1]^k : \sum_{j=1}^k x_j \les  1 \biggr\}.
	$$ We put 
	$\omega_n = 0$ for all $n$ except $n \equiv \nu_0 \pmod{W}$ for some fixed 
	$\nu_0$ such that $(\nu_0 + h_j, W) = 1$ for all $j$. We also put 
	$\lambda_{d_1, \ldots, d_k} = 0$ if 
	$\bigr( \prod_{j=1}^k d_j, W \bigr) > 1$ for at least one $j$.  
	
	We obtain the desired result using the following assertion 
	(see \cite{Maynard}):
	
	1) Under the above assumptions on $\omega_n$, the following 
	relations hold:
	\begin{gather*}
	S_1 = \sum_{\substack{X \les  n < 2X \\ n \equiv \nu_0 \pmod{W}}} 
	\omega_n = \frac{\bigl( 1 + o(1) \bigr) \varphi^k (W) X 
		(\log R)^k}{W^{k+1}} I_k (F), \\
	S_2 = \sum_{\substack{X \les  n < 2X \\ n \equiv \nu_0 \pmod{W}}} 
	\biggl( \sum_{j=1}^k \chi_{\mathbb{P}} (n+h_j) \biggr) \omega_n = \frac{
		\bigl( 1 + o(1) \bigr) \varphi^k (W) X (\log R)^{k+1}}{W^{k+1} \log X} 
	\sum_{i=1}^k J_k^{(j)} (F),
	\end{gather*} provided $I_k (F) \ne 0$ and $J_k^{(j)}(F) \ne 0$ for each $j$, where
	\begin{gather*}
	I_k (F) = \int_0^1 \ldots \int_0^1 F(t_1, \ldots, t_k)^2 dt_1 \ldots dt_k, \\
	J_k^{(j)}(F) = \int_0^1 \ldots \int_0^1 \biggl( \int_0^1 F(t_1, \ldots, 
	t_k) dt_j \biggr)^2 dt_1 \ldots dt_{j-1} dt_{j+1} \ldots dt_k.
	\end{gather*}
	
	2) Define
	$$
	M_k = \sup_{F} \frac{\sum_{j=1}^k J_k^{(j)} (F)}{I_k (F)}.
	$$ Then for all $k \ges 600$ the following inequality holds true:
	\begin{equation}\label{M_k}
	M_k > \log k - 2 \log \log k - 1.
	\end{equation}
	
	To apply Maynard's argument we need an analogue of part 1)
	for $S_{1,\alpha}$ and $S_{2,\alpha}$. 
	It would follow from the relations
	\begin{equation} \label{relation_to_Maynard}
	S_{1,\alpha} = \biggl( d - c + o(1) \biggr) S_1, \qquad S_{2,\alpha} = 
	\biggl( d - c + o(1) \biggr) S_2.
	\end{equation} Here the coefficient $d - c < 1$ corresponds to the density of 
	$\mathbb{E}$ among the integers. Note that since $\alpha < 1$, the numbers $\{ n^{\alpha} 
	\}, \{ (n+h_1)^{\alpha} \}, \ldots, \{ (n+h_k)^{\alpha} \}$ are close 
	to each other. So, in order to verify that all $k+1$ conditions $n \in 
	\mathbb{E}$, 
	$n+h_j \in \mathbb{E}$, $1 \les  j \les  k$ hold true, it is 
	sufficient to check only 
	two of them: $n \in \mathbb{E}, n+h_k \in \mathbb{E}$. Thus, we define
	the new subset as
	$$
	\mathbb{E}' = \bigl\{ n \in \mathbb{N} : \{ n \in \mathbb{E} \} \cap 
	\{ n+h_k \in \mathbb{E} \} \bigr\}.
	$$ Then
	$$
	S_{2,\alpha} = \sum_{j=1}^k S_{2,\alpha}^{(j)} = \sum_{j=1}^k 
	\sum_{\substack{X \les  n < 
			2X \\ n \in \mathbb{E}'}} \chi_{\mathbb{P}} (n+h_j) \omega_n.
	$$ Following~\cite{Maynard}, we change the order of summation in $S_{2,
		\alpha}^{(j)}$ and apply Chinese Remainder theorem. Thus we get
	$$
	S_{2,\alpha}^{(j)} = \sum_{\substack{d_1, \ldots, d_k \\ e_1, \ldots, 
			e_k \\ ([d_i, e_i], [d_j, e_j]) = 1 \ \forall i \ne j}} 
	\lambda_{d_1, \ldots, d_k} \lambda_{e_1, \ldots, e_k} 
	\sum_{\substack{X \les  n < 2X \\ n \equiv a \pmod{q} \\ n \in 
			\mathbb{E}'}} 
	\chi_{\mathbb{P}}(n + h_j),
	$$ where $q = W \prod_{i=1}^k [d_i, e_i]$. Further,
	\begin{multline} \label{S2alpha}
	S_{2,\alpha}^{(j)} = \frac{\pi_{\mathbb{E}} (2X) - \pi_{\mathbb{E}} (X) + 
		O(1)}{\varphi(W)} 
	\sum_{\substack{d_1, \ldots, d_k \\ e_1, \ldots, e_k \\ 
			e_j = d_j = 1 \\ ([d_i, e_i], [d_j, e_j]) = 1 \ \forall i \ne j}}
	\frac{\lambda_{d_1, \ldots, d_k} \lambda_{e_1, 
			\ldots, e_k}}{\prod_{i=1}^k \varphi \biggl( [d_i, e_i] \biggr)} 
	+ \\ O \biggl( \sum_{
		\substack{d_1, \ldots, d_k \\ e_1, \ldots, e_k}} |\lambda_{d_1, 
		\ldots, d_k} \lambda_{e_1, \ldots, e_k}| E^{(j)} (X, q) \biggr),
	\end{multline} where $\pi_{\mathbb{E}}$ is the counting
	function of primes from $\mathbb{E}$,
	\begin{align*}
	E^{(j)} (X,q) = 1 + \max_{(a,q)=1} \biggl| \sum_{\substack{X \les  
			n < 2X \\ 
			n \equiv a\pmod{q} \\ n \in \mathbb{E}'}}
	\chi_{\mathbb{P}}(n+h_j) - \frac{1}{\varphi(q)} &\sum_{\substack{X 
			\les  n < 2X \\ n \in \mathbb{E}'}} \chi_{\mathbb{P}} (n+h_j) 
	\biggr| \les  \\  1 + E_1^{(j)} (&X, q) + E_2^{(j)} (X,q) 
	+ E_3^{(j)} (X,q), \\
	E_1^{(j)} (X, q) = \max_{(a,q)=1} \sum_{\substack{X \les  
			n < 2X \\ n 
			\equiv a \pmod{q} \\ n \in \mathbb{E} \setminus \mathbb{E}'}} 
	\chi_{\mathbb{P}}(n+h_j), \qquad
	E_2^{(j)} (X&,q) = \frac{1}{\varphi(q)} 
	\sum_{\substack{X \les  n < 2X \\ n \in \mathbb{E} \setminus 
			\mathbb{E}'}} 
	\chi_{\mathbb{P}} (n + h_j), \\
	E_3^{(j)} (X, q) = \max_{(a,q)=1} \biggl| \sum_{\substack{X \les  n < 2X \\ 
			n \equiv a \pmod{q} \\ n \in \mathbb{E}}} \chi_{\mathbb{P}}(n+h_j) 
	- &\frac{1}{\varphi(q)} \sum_{\substack{X \les  n < 2X \\ n \in 
			\mathbb{E}}} 
	\chi_{\mathbb{P}}(n+h_j) \biggr|.
	\end{align*} The trivial upper bound for $E_1^{(j)}$, $E_2^{(j)}$ and $E_3^{(j)}$ is 
	$X / \varphi(q)$. Similarly 
	to (5.20) from~\cite{Maynard}, we apply Cauchy's inequality to the error term in~\eqref{S2alpha} to bound it from above by
	$$
	\lambda_{\max}^2 \biggl( \sum_{q < R^2 W} \mu^2 (q) \tau_{3k}^2 (q) \frac{X}{
		\varphi(q)} \biggr)^{1/2} \biggl( \sum_{q < R^2 W}   \bigl(1 + E_1^{(j)} 
	(X, q) + E_2^{(j)} (X, q) + E_3^{(j)} (X,q) \bigr) \biggr)^{1/2},
	$$ where $\lambda_{\max} = \max_{d_1, \ldots, d_k} 
	|\lambda_{d_1, \ldots, d_k}|$.
	The sum in the second factor doesn't exceed $X (\log X)^{-B}$ for any 
	fixed of $B > 0$. Indeed, \hyperref[thm2]{Theorem~2} implies that the contribution coming
	from the term $E_3^{(j)}(X,q)$ is estimated by $X(\log X)^{-A}$ 
	for any fixed $A > 0$. 
	We estimate $E_1^{(j)} (X,q)$ and $E_2^{(j)} (X,q)$ by the way similar to the estimate
	of $S^{(2)}$ and $S^{(3)}$ in \hyperref[sec2]{Section 2} 
	(note that all the primes $n+h_j$ which belong 
	to $\mathbb{E}$ but do not belong to $\mathbb{E}'$ lie in the subset 
	$\mathbb{E}_{\Delta}$ defined earlier). The sum in the first factor doesn't 
	exceed $X (\log X)^{c_k}$ for some $c_k > 0$. So, choosing $A$ large enough, 
	we get the bound 
	of $\lambda_{\max}^2 X (\log X)^{-B}$ for the error term in~\eqref{S2alpha}.
	Finally, we apply Vinogradov's result~\cite{Vinogradov_1} to 
	$\pi_{\mathbb{E}} (X)$ to get the second relation 
	in~\eqref{relation_to_Maynard}.
	
	We treat the sum $S_{1, \alpha}$ in a similar way. Thus we get
	\begin{multline} \label{S_1_asymp}
	S_{1,\alpha} = \sum_{\substack{d_1, \ldots, d_k \\ e_1, \ldots, 
			e_k \\ ([d_i, e_i], [d_j, e_j]) = 1 \ \forall i \ne j}} \lambda_{d_1, 
		\ldots, d_k} \lambda_{e_1, \ldots, e_k} \sum_{\substack{X \les  n < 
			2X \\ n \equiv a \pmod{q} \\ n \in \mathbb{E}'}} 1 = \\
	\frac{X}{
		2W} \sum_{\substack{d_1, \ldots, d_k \\ e_1, \ldots, e_k \\
			([d_i, e_i], [d_j, e_j]) = 1 \ \forall i \ne j}} \frac{
		\lambda_{d_1, \ldots, d_k} \lambda_{e_1, \ldots, e_k}}{\prod_{i=1}^k 
		[d_i, e_i]} + O \biggl( \sum_{\substack{d_1, \ldots, 
			d_k \\ e_1, \ldots, e_k}} \lambda_{\max}^2 G(X,q) \biggr),
	\end{multline} where
	\begin{gather*}
	G(X,q) = 1 + \max_{(a,q)=1} \biggl| \sum_{\substack{X \les  n < 2X \\ n 
			\equiv 
			a\pmod{q} \\ n \in \mathbb{E}'}} 1 - \frac{1}{q} \sum_{\substack{X 
			\les  n < 2X \\ n \in \mathbb{E}'}} 1 \biggr| \les  1 + G_1 
	(X,q) + G_2 (X,q) + G_3 (X,q), \\
	G_1 (X,q) = \max_{(a,q)=1} \sum_{\substack{X \les  n < 2X \\ 
			n \equiv a \pmod{q} \\ n \in \mathbb{E} \setminus \mathbb{E}'}} 1 
	, \qquad
	G_2 (X,q) = \frac{1}{q} \sum_{\substack{X 
			\les  n < 2X \\ n \in \mathbb{E} \setminus \mathbb{E}'}} 1, \\
	G_3 (X,q) = \max_{(a,q) = 1} \biggl| \sum_{\substack{X \les  n < 
			2X \\ n \equiv a \pmod{q} \\ n \in \mathbb{E}}} - \frac{1}{q} 
	\sum_{\substack{X \les  n < 2X \\ n \in \mathbb{E}}} 1 \biggr|.
	\end{gather*} We estimate the error term in~\eqref{S_1_asymp} similarly to the error term
	in~\eqref{S2alpha}. For $G_3$ we apply the arguments used in the proof of
	\hyperref[thm2]{Theorem 2}. This case is simpler since the summation goes over integers so 
	one deals with
	$$
	\sum_{\substack{X \les  n < 2X \\ n \equiv a \pmod{q}}} e \bigl(h n^{\alpha} 
	\bigr).
	$$ The congruence condition is removed by substitution $n = qr + l$. 
	The sum over $r$ is then estimated by~\eqref{Heath-Brown} with $k=3$. 
	This concludes the proof of~\eqref{relation_to_Maynard}.
	We note that the estimates for $E_3^{(j)}$ and $G_3$ 
	correspond to the conditions (2) and (1) in the Hypothesis 1 
	from~\cite{Maynard_2}.
	
	We show that there are infinitely many $n$ such that at least 
	$\lceil M_k / 6 \rceil$ numbers $n+h_i$ are primes 
	from $\mathbb{E}$. By definition of $M_k$, there is a function $F_0$ such that
	$$
	\sum_{j=1}^k J_k^{(j)}(F_0) > (M_k - \delta_1) I_k (F_0).
	$$ Then
	$$
	S_{\alpha} > (d-c) \frac{\varphi^k (W) X 
		(\log R)^k I_k (F_0)}{W^{k+1}} \biggl( \biggl( \frac{1}{6} - 
	\delta_1 \biggr) (M_k - \delta_1) - \rho + o(1) \biggr).
	$$ If $\rho = M_k / 6 - \delta_2$ for some $\delta_2 > 0$
	such that $\lceil M_k / 6 \rceil = \lfloor \rho + 1 \rfloor$ 
	and $\delta_1$ is small
	enough (depending on $\delta_2$), then there are infinitely many $n$ 
	such that at least $\lfloor 
	\rho + 1 \rfloor = m+1$ numbers among $n+h_j$ are primes. 
	
	To finish the computation, we estimate the right hand side of~\eqref{maximum}. 
	In accordance with~\eqref{M_k}, we choose $k = \lfloor 390 m^2 
	e^{6m} \rfloor$ and take the tuple of primes 
	$\{p_{\pi(k) + 1}, \ldots, p_{\pi(k)+k} \}$ 
	which is obviously admissible (indeed, there are no element of this set
	which is congruent to 
	zero modulo any prime $p \les  k$; on the other hand, is does not cover a 
	complete residue system modulo any $p > k$ due to its size). Next,
	$$
	\max_{1 \les  i < j \les  k} |h_j - h_i|
	\les  p_{\pi(k) + k} \les  p_{\lceil 1.1 k \rceil}
	$$ for $k \ges 10^5$. Using the inequality $p_n < n (\log n + 
	\log \log n + 8)$ 
	(see, for example,~\cite{Rosser}) with $n=\lceil 1.1\cdot 390 
	m^2 e^{6m} \rceil$, we find that
	$$
	\liminf_{n \to +\infty} (p_{n+m} - p_n) \les  9~700 m^3 e^{6m}.
	$$ This bound can certainly be sharpened for small $m$. For example, 
	if $m=1$ then we take $k = 157~337$. Thus we get $M_k > 6$,
	and the precise computation of the right hand side of~\eqref{maximum} gives
	$$
	p_{\pi(k) + k} - p_{\pi(k)+1} = p_{171~807} - p_{14~471} = 2~176~652 
	< 2.18 \cdot 10^6.
	$$ For $m=2$, we need $M_k > 12$, so we take $k = 157~629~323$ and
	$$
		p_{\pi(k) + k} - p_{\pi(k) + 1} = p_{166~478~324} - p_{8~849~002} = 
		3~130~607~572 < 3.14 \cdot 10^9.
	$$

	\section*{ACKNOWLEDGMENTS}
	
	I would like to thank M.A.~Korolev for substantial help throughout the 
	preparation of this paper. I also thank M.~Radziwill for fruitful 
	discussions.

\end{document}